\documentclass[12pt]{amsart}

\usepackage{amsmath}
\usepackage{amssymb}
\usepackage{array}

\renewcommand{\atop}[2]{%
\genfrac{}{}{0pt}{}{#1}{#2}}

\newtheorem{theorem}{Theorem}
\newtheorem{lemma}{Lemma}
\newtheorem{corollary}{Corollary}

\theoremstyle{definition}

\renewcommand{\atop}[2]{%
\genfrac{}{}{0pt}{}{#1}{#2}}

\begin{document}

\title{Approximations to Euler's constant}

\author{Kh.~Hessami Pilehrood$^{1}$}

\address{Institute for Studies in Theoretical Physics and Mathematics
(IPM), Tehran, Iran} \curraddr{Mathemetics Department, Faculty of
Science, Shahrekord University, Shahrekord, P.O. Box 115, Iran.}
\email{hessamik@ipm.ir, hessamit@ipm.ir, hessamit@gmail.com}
\thanks{$^1$ This research was in part supported by a grant
from IPM (No. 86110025)}

\author{T.~Hessami Pilehrood$^2$}
\thanks{$^2$ This research was in part supported by a grant
from IPM (No. 86110020)}

\subjclass{11Y60, 11Y35, 41A25}

\date{}

\keywords{}

\begin{abstract}
We study a problem of finding good approximations to Euler's
constant $\gamma=\lim_{n\to\infty}S_n,$ where
$S_n=\sum_{k=1}^n\frac{1}{n}-\log(n+1),$ by linear forms in
logarithms and harmonic numbers. In 1995, C.~Elsner showed that
 slow convergence of the sequence $S_n$ can be significantly
improved if  $S_n$ is replaced by linear combinations of $S_n$
with integer coefficients. In this paper, considering more general
linear transformations of the sequence $S_n$ we establish new
accelerating convergence formulae for $\gamma.$ Our estimates
sharpen and generalize recent Elsner's, Rivoal's and author's
results.
\end{abstract}

\maketitle

\section{Introduction}
\label{intro} Let $\alpha\ge 0$ be a real number and
$$\gamma_{\alpha}=\sum_{k=1}^{\infty}\left(\frac{1}{k+\alpha}-\log\left(\frac{k+\alpha+1}{k+\alpha}\right)
\right).
$$
We denote the partial sum of the above series by
\begin{equation}
\begin{split}
S_n(\alpha)&=\sum_{k=1}^n\left(\frac{1}{k+\alpha}-\log\left(\frac{k+\alpha+1}{k+\alpha}\right)
\right) \\ &=\sum_{k=1}^n\frac{1}{k+\alpha} -\log(\alpha+n+1)
+\log(\alpha+1) \label{eq1}
\end{split}
\end{equation}
and $S_n:=S_n(0).$ It easily follows (see \cite[formula (2)]{ta})
that
$$
\lim_{n\to\infty}S_n(\alpha)=-\frac{\Gamma'(\alpha+1)}{\Gamma(\alpha+1)}+\log(\alpha+1)=
-\psi(\alpha+1)+\log(\alpha+1),
$$
where $\psi(\alpha)$ is the logarithmic derivative of the gamma
function (or the digamma function) and therefore,
$$
\gamma_{\alpha}=\log(\alpha+1)-\psi(\alpha+1).
$$
In particular, $\gamma_0=-\psi(1)=\gamma=0.577215\ldots,$ where
$\gamma$ is Euler's constant. It is well-known that the sequence
$S_n$ slowly converges to the Euler constant $\gamma$ (see,  for
details, \cite{ka})
$$
\gamma=S_n+O(n^{-1}).
$$
In 1995, Elsner \cite{el} found out that $\gamma$ can be
approximated by linear combinations of partial sums (\ref{eq1})
with integer coefficients
\begin{equation}
\left|\gamma-\sum_{k=0}^n(-1)^{n+k}\binom{n}{k}
\binom{k+n+\tau-1}{k+\tau-1}S_{k+\tau-1}\right|\le\frac{1}%
{2n\tau\binom{n+\tau}{n}}, \qquad \tau,n \in\mathbb{N} \label{eq2}
\end{equation}
and this inequality exhibits geometric convergence if $\tau=O(n).$
Formulas (\ref{eq2}) for $\tau>n$ were generalized by Rivoal in
\cite{ri}, where, in particular, it was shown that
$$
\left|\gamma-\frac{1}{2^n}\sum_{k=0}^n(-1)^{k+n}\binom{n}{k}\binom{2k+2n}{n}
S_{2k+n}\right|=O\left(\frac{1}{n27^{n/2}}\right), \quad
n\to\infty.
$$
Another such kind formula
$$
\gamma-\sum_{k=0}^n(-1)^{k+n}\binom{n}{k}\binom{n+k}{k}S_{k+n}=\frac{1}{4^{n+o(n)}},
\quad n\to\infty
$$
was proved in \cite{he}. Recently, C.~Elsner \cite{el2} presented
a two-parametric series transformation of the sequence $S_n$
\begin{equation}
\sum_{k=0}^n(-1)^{n+k}\binom{n}{k}\binom{n+\tau_1+k}{n}S_{k+\tau_2-1}
\label{eq3}
\end{equation}
converging more rapidly to $\gamma$ when $\tau_2>\tau_1+1$ and $n$
increases than in the case $\tau_2=\tau_1+1$ considered in
(\ref{eq2}).

In this paper, we consider a more general series transformation of
the type
\begin{equation}
\frac{n_1!\ldots n_m!}{N!\, r^N}\sum_{k=0}^N(-1)^{N+k}\binom{N}{k}
\binom{rk+n_1+\tau_1}{n_1}\cdots \binom{rk+n_m+\tau_m}{n_m}
S_{rk+\tau_0} \label{3.5}
\end{equation}
with $n_1,\ldots, n_m\in {\mathbb N},$ $\tau_0, \tau_1,\ldots,
\tau_m\in {\mathbb N}_0,$ and $N=\sum_{j=1}^mn_j,$ and
 give new accelerating convergence formulae for Euler's
constant $\gamma.$  In particular,  we show  (see Theorem \ref{t2}
and Corollary \ref{c1} below) that  if $\tau_1, \tau_2$ are linear
functions of $n,$ then the sum (\ref{eq3}) converges to $\gamma$
at the least geometric rate and represents the best approximation
in the set of all the sums (\ref{eq3}) with a fixed value of
$\lim\limits_{n\to\infty}\tau_2/n,$ provided that
$\lim\limits_{n\to\infty}2(\tau_2-\tau_1)/n=1.$

\section{Statement of the main results}
\label{st}

As usual, we denote the Gauss  hypergeometric function (see, for
details, \cite{ra}) by
$$
{}\sb 2F\sb
 {1}\left(\left.\atop{a,
b}{c}\right|z\right)=\sum_{\nu=0}^{\infty} \frac{(a)_{\nu}
(b)_{\nu}}{\nu\,!\,(c)_{\nu}}{z^{\nu}},
$$
 where $(\lambda)_{\nu}$ is the Pochhammer symbol (or the shifted
factorial) defined by
\begin{equation*}
(\lambda)_{\nu}=\frac{\Gamma(\lambda+\nu)}{\Gamma(\lambda)}
=\begin{cases}
 1,     & \quad \nu=0; \\
\lambda(\lambda+1)\ldots (\lambda+\nu-1), & \quad \nu\in {\mathbb
N}.
\end{cases}
\end{equation*}
We then prove the following theorems:
\begin{theorem} \label{t1}
Let $n_1,\ldots, n_m\in {\mathbb N},$ $\tau_0, \tau_1,\ldots,
\tau_m\in {\mathbb N}_0,$
 $0\le\tau_0-\tau_m\le n_m,$
$n_m+\tau_m\ge n_j+\tau_j,$ $j=1,\ldots, m-1,$ and
$N=\sum_{j=1}^mn_j.$
  Then
\begin{equation}
\begin{split}
 &\left|\frac{N!\,(-r)^N}{n_1!\ldots n_m!}\,\gamma-\sum_{k=0}^N(-1)^{k}\binom{N}{k}
\binom{rk+n_1+\tau_1}{n_1}\cdots \binom{rk+n_m+\tau_m}{n_m}
S_{rk+\tau_0}\right|
\\ &=\prod_{j=1}^{m}
\binom{n_m+\tau_m-\tau_j}{n_j}\int\limits_0^1\int\limits_0^1
\frac{x^{n_m+\tau_m}(1-x^r)^Nt^{n_m+\tau_m-\tau_0}(1-t)^{\tau_0-\tau_m}\omega(t)}%
{(1-t+xt)^{n_m+1}} \\
&\times \left|Q_m\left(\frac{xt}{1-t+xt}\right)\right|\, dxdt,
\label{eq4}
\end{split}
\end{equation}
where
\begin{equation}
\omega(t)=\frac{1}{t(\log^2(1/t-1)+\pi^2)} \label{eq5}
\end{equation}
and $Q_m(y)$ is a polynomial of degree $N-n_m$ given by the
formula
\begin{equation}
Q_m(y)=\sum_{k_1=0}^{n_1}\ldots\sum_{k_{m-1}=0}^{n_{m-1}}
\prod_{j=1}^{m-1}\frac{(-n_j)_{k_j}(1+n_m+\tau_m-\tau_{j+1})_{k_1+\ldots+k_j}}%
{\,\,\,k_j\,!\,(1+n_m+\tau_m-n_j-\tau_j)_{k_1+\ldots+k_j}}y^{k_j}
\label{q}
\end{equation}
if $m\ge 2,$ and $Q_1(y)\equiv 1.$
\end{theorem}
\begin{theorem} \label{t2}
Let $b,c,r\in {\mathbb N},$ $a\in {\mathbb N}_0,$ $0\le b-a\le c.$
 Then for $n\in {\mathbb N}$ we have
\begin{equation}
\left|\gamma-\frac{1}{r^{cn}}
\sum_{k=0}^{cn}(-1)^{k+cn}\binom{cn}{k}
\binom{rk+(a+c)n}{cn}S_{rk+bn}\right|
<\left(\frac{b^{\frac{b}{r}}(c+a-b)^{c+a-b}(b-a)^{b-a}}{(b+cr)^{c+
\frac{b}{r}}}\right)^n
 \label{tt2}
\end{equation}
 (Here and throughout the paper $0^0$ is treated as
$1.$)
\end{theorem}
 If $b,c,r$  are fixed, then the minimum of the right-hand side
of (\ref{tt2}) is attained when $b-a=c/2$ and in this case we have
\begin{corollary} \label{c1}
Let $b,c,r,n\in {\mathbb N}$ and $b\ge c.$  Then
\begin{equation*}
\left|\gamma-\frac{1}{r^{2cn}}
\sum_{k=0}^{2cn}(-1)^{k}\binom{2cn}{k}
\binom{rk+(b+c)n}{2cn}S_{rk+bn}\right|
<\left(\frac{b^{\frac{b}{r}}c^{2c}}{(b+2cr)^{2c+
\frac{b}{r}}}\right)^n.
\end{equation*}
\end{corollary}
\begin{theorem} \label{t3}
Let $b,c,r\in {\mathbb N},$ $a\in {\mathbb N}_0$ and $0\le b-a\le
c.$ Then for any positive integer $n\ge 2/c$ one has
\begin{equation*}
\begin{split}
&\left|\gamma-\frac{((cn)!)^2}{(2cn)!\,r^{2cn}}
\sum_{k=0}^{2cn}(-1)^{k}\binom{2cn}{k}
\binom{rk+(a+c)n}{cn}^2S_{rk+bn}\right|\\
&<cn\left(\frac{b^{\frac{b}{r}}c^c(c+a-b)^{c+a-b}(b-a)^{b-a}}{(b+2cr)^{2c+
\frac{b}{r}}}\right)^n.
\end{split}
\end{equation*}
\end{theorem}
By the similar argument as above putting $a=b-c/2$ we get a
sharper bound than in Corollary \ref{c1}.
\begin{corollary} \label{c2}
Let $b,c,r,n\in {\mathbb N},$ $2b\ge c,$ and $c$ is even.  Then
\begin{equation*}
\left|\gamma-\frac{((cn)!)^2}{(2cn)!\,r^{2cn}}
\sum_{k=0}^{2cn}(-1)^{k}\binom{2cn}{k}
\binom{rk+\left(b+\frac{c}{2}\right)n}{cn}^2S_{rk+bn}\right|
<cn\left(\frac{b^{\frac{b}{r}}c^{2c}}{2^c(b+2cr)^{2c+
\frac{b}{r}}}\right)^n.
\end{equation*}
\end{corollary}
For example, setting $b=c=4, r=1$ we get the following estimate:
\begin{corollary} \label{c2.5}
For any positive integer $n$ one has
\begin{equation*}
\left|\gamma-\frac{(4n)!^2}{(8n)!}\sum_{k=0}^{8n}(-1)^k\binom{8n}{k}
\binom{k+6n}{4n}^2S_{k+4n}\right|<\frac{4n}{(2^4\ 3^{12})^n}<4n
(0.00000012)^n.
\end{equation*}
\end{corollary}
\begin{theorem} \label{t4}
Let $n_1,\ldots, n_m\in {\mathbb N},$ $\tau_0, \tau_1,\ldots,
\tau_m\in {\mathbb N}_0,$
 $0\le\tau_0-\tau_m\le n_m,$
$n_m+\tau_m\ge\tau_{j+1}> n_j+\tau_j,$ $j=1,\ldots, m-1,$ and
$N=\sum_{j=1}^mn_j,$
  Then
\begin{equation*}
\begin{split}
 &\left|\frac{N!\,(-r)^N}{n_1!\ldots n_m!}\,\gamma-\sum_{k=0}^N(-1)^{k}\binom{N}{k}
\binom{rk+n_1+\tau_1}{n_1}\cdots \binom{rk+n_m+\tau_m}{n_m}
S_{rk+\tau_0}\right|
\\ &\le\prod_{j=1}^{m}\!\!
\binom{n_m+\tau_m-\tau_j}{n_j}\!\!\int\limits_0^1\!\!\!\int\limits_0^1\!
\frac{x^{n_m+\tau_m}(1-x^r)^Nt^{n_m+\tau_m-\tau_0}(1-t)^{\tau_0-\tau_m}\omega(t)}%
{(1-t+xt)^{n_m+1}}\, dxdt.
\end{split}
\end{equation*}
\end{theorem}
Setting $\tau_{j+1}=n_j+\tau_j+1,$ $j=1,\ldots,m-1,$ in Theorem
\ref{t4} we get
\begin{corollary} \label{c3}
Let $n_1, \ldots, n_m\in {\mathbb N},$ $\tau_0,\tau_1\in {\mathbb
N}_0,$ $N=\sum_{j=1}^mn_j,$ and $N-n_m+\tau_1+(m-1)\le\tau_0\le
N+\tau_1+(m-1).$ Then
\begin{equation*}
\begin{split}
&\left|\gamma-\frac{n_1!\ldots
n_m!}{N!(-r)^N}\sum_{k=0}^N(-1)^k\binom{N}{k}\prod_{j=1}^m
\binom{\!rk+n_1+\ldots+n_j+\tau_1+j-1\!}{n_j} S_{rk+\tau_0}\right|\\
&\le\prod_{j=1}^{m-1}\frac{N+j}{n_{j+1}+\ldots+n_m+m-j}\times\\
&\int\limits_0^1\int\limits_0^1
\frac{x^{N+\tau_1+m-1}(1-x^r)^Nt^{N+\tau_1+m-1-\tau_0}
(1-t)^{\tau_0+n_m-N-\tau_1-m+1}\omega(t)}%
{r^N(1-t+xt)^{n_m+1}}\,dxdt
\end{split}
\end{equation*}
\end{corollary}
\begin{theorem} \label{t5}
Let $m,c_1,\ldots,c_m,r,b,n\in {\mathbb N},$ $a\in {\mathbb N}_0,$
$C=\sum_{j=1}^mc_j,$ and $a-c_m\le b-c\le a.$ Then
\begin{equation*}
\begin{split}
&\left|\gamma-\frac{(c_1n)!\ldots (c_mn)!}{(Cn)!(-r)^{Cn}}
\sum_{k=0}^N(-1)^k\binom{Cn}{k}\prod_{j=1}^m
\binom{\!rk+(a+c_1+\ldots+c_j)n+j-1\!}{c_jn}\right.\\
&\times\Biggl. S_{rk+bn+m}\Biggr|
<M(\overline{c})\left(\frac{b^{\frac{b}{r}}C^C(C+a-b)^{C+a-b}
(c_m+b-a-C)^{c_m+b-a-C}}{c_m^{c_m}(b+Cr)^{C+\frac{b}{r}}}\right)^n,
\end{split}
\end{equation*}
where $M(\overline{c})<C^{m-1}$ is some constant depending only on
$c_1, \dots, c_m.$
\end{theorem}
Consider several illustrative examples of Theorem \ref{t5}. Taking
$c_1=\ldots=c_m=2c,$ $C=2mc,$ $b=2mc,$ $a=c,$ $c\in {\mathbb N},$
we get
\begin{corollary} \label{c4}
Let $c,m,r\in {\mathbb N}.$ Then for any positive integer $n$ one
has
\begin{equation*}
\begin{split}
&\left|\gamma-\frac{(2cn)!^m}{(2mcn)!r^{2mcn}}\sum_{k=0}^{2mcn}(-1)^k
\binom{2mcn}{k}\binom{rk+3cn+1}{2cn}\binom{rk+5cn+2}{2cn}\ldots\right.\\
&\left.\times\binom{rk+(2m+1)cn+m}{2cn}S_{rk+2mcn+m}\right|
<\frac{m^m}{(m-1)!}\left(\frac{1}{4^c(r+1)^{2mc+\frac{2mc}{r}}}\right)^n
\end{split}
\end{equation*}
\end{corollary}
Setting $c_1=\ldots=c_m=2c,$ $C=2mc,$ $b=(2m-1)c,$ $a=2c,$ $c\in
{\mathbb N},$  we get
\begin{corollary}
Let $c,m,r\in {\mathbb N}.$ Then for any positive integer $n$ one
has
\begin{equation*}
\begin{split}
&\left|\gamma-\frac{(2cn)!^m}{(2mcn)!r^{2mcn}}\sum_{k=0}^{2mcn}(-1)^k
\binom{2mcn}{k}\prod_{j=1}^m\binom{rk+2jcn+j}{2cn}
S_{rk+(2m-1)cn+m}\right| \\
&<\frac{m^m}{(m-1)!}
\left(\frac{4^{-c}\left(1-\frac{1}{2m}\right)^{\frac{(2m-1)c}{r}}}{\left(
r+1-\frac{1}{2m}\right)^{2mc+\frac{(2m-1)c}{r}}}\right)^n.
\end{split}
\end{equation*}
\end{corollary}

\section{Analytical construction}

We define the generalized Legendre polynomial by
$A(x)=\sum_{k=0}^NA_kx^{rk}$ with
$$
A_k=(-1)^{k+N}\binom{N}{k}\binom{rk+n_1+\tau_1}{n_1}\cdots\binom{rk+n_m+\tau_m}{n_m}.
$$
\begin{lemma} \label{l1}
There holds
$$
A(1)=\sum_{k=0}^NA_k=\frac{N!\,r^N}{n_1!\ldots n_m!}.
$$
\end{lemma}

{\bf Proof.} For the proof, let
$$
R(t)=\frac{N!}{n_1!\ldots
n_m!}\frac{(rt-n_1-\tau_1)_{n_1}(rt-n_2-\tau_2)_{n_2}\ldots
(rt-n_m-\tau_m)_{n_m}}{t(t+1)\ldots (t+N)}.
$$
Such rational functions were considered early by the authors
\cite{he1}, \cite{he2} to derive explicit Pad\'e approximations of
the first and second kinds for polylogarithmic functions. As it is
easily seen the rational function $R(t)$ has the following
partial-fraction expansion:
$$
R(t)=\sum_{k=0}^N\frac{A_k}{t+k},
$$
from which it follows that
$$
\sum_{k=0}^NA_k=\sum_{k=0}^N\underset{t=-k}{{\rm res}} R(t)=
-\underset{t=\infty}{{\rm res}} R(t)=\frac{N!\,r^N}{n_1!\ldots
n_m!}.  \qquad \qquad \qed
$$
Put
$$
I(\alpha):=\int_0^1x^{\tau_0+\alpha}A(x)\left(\frac{1}{1-x}+\frac{1}{\log
x}\right)\,dx
$$
\begin{lemma} \label{l2}
There holds the equality
$$
I(\alpha)=\frac{N!\,r^N}{n_1!\ldots
n_m!}\,\gamma_{\alpha}-\sum_{k=0}^NA_kS_{rk+\tau_0}(\alpha).
$$
\end{lemma}
{\bf Proof.} Substituting
$$
\frac{1}{1-x}+\frac{1}{\log x}=\int_0^1\frac{1-x^t}{1-x}\,dt,
$$
we get
$$
I(\alpha)=\int_0^1\int_0^1x^{\tau_0+\alpha}A(x)\frac{1-x^t}{1-x}\,dtdx=
\sum_{k=0}^NA_k\int_0^1\int_0^1\frac{x^{rk+\tau_0+\alpha}(1-x^t)}{1-x}\,dxdt.
$$
Expanding $(1-x)^{-1}$ in a geometric series and applying Lemma
\ref{l1} we find
\begin{equation*}
\begin{split}
I(\alpha)&=\sum_{k=0}^NA_k\sum_{l=0}^{\infty}\int_0^1\int_0^1x^{rk+\tau_0+l+\alpha}
(1-x^t)\,dxdt \\
&=\sum_{k=0}^NA_k\sum_{l=0}^{\infty}\int_0^1\left(\frac{1}{rk+\tau_0+l+\alpha+1}-\frac{1}%
{rk+\tau_0+t+l+\alpha+1}\right)\,dt \\
&=\sum_{k=0}^NA_k\sum_{l=1}^{\infty}\left(\frac{1}{rk+\tau_0+l+\alpha}-\log\left(
\frac{rk+\tau_0+l+\alpha+1}{rk+\tau_0+l+\alpha}\right)\right) \\
&=\sum_{k=0}^NA_k(\gamma_{\alpha}-S_{rk+\tau_0}(\alpha))=
\frac{N!\,r^N}{n_1!\ldots
n_m!}\,\gamma_{\alpha}-\sum_{k=0}^NA_kS_{rk+\tau_0}(\alpha).
\qquad\qquad\qed
\end{split}
\end{equation*}
Next, we consider two differential operators
\begin{equation*}
\begin{split}
S_{\tau,n}(f(x))&=\frac{(-1)^n}{n!}x^{-\tau}\left(x^{n+\tau}f(x)\right)^{(n)},
\\
T_{\tau,n}(f(x))&=\frac{1}{n!}x^{n+\tau}\left(x^{-\tau}f(x)\right)^{(n)},
\end{split}
\end{equation*}
where $\tau$ is a real number and $n$ is a non-negative integer.
We show that $S_{\tau,n}$ and $T_{\tau,n}$ are adjoint operators
in  some sense.
\begin{lemma} \label{l3}
Suppose that $f(x)$ is a polynomial vanishing at $x=1$ with order
at least $n$ and $g(x)\in C^{\infty}(0,1)\cap L^1(0,1)$ satisfies
the following boundary conditions:
$$
\lim_{x\to 0+}x^lg^{(l-1)}(x)=\lim_{x\to 1-}(1-x)^lg^{(l-1)}(x)=0
$$
for all $1\le l\le n.$ Then we have
$$
\int_0^1S_{\tau,n}(f(x))\cdot g(x)\,dx=\int_0^1f(x)\cdot
T_{\tau,n}(g(x))\,dx.
$$
\end{lemma}

{\bf Proof.}  The proof is analogous to the proof of Lemma 3.1
\cite{ha}. \qed
\begin{lemma} \label{l4}
There holds
$$
I(\alpha)=\int\limits_0^1\!\!\int\limits_0^1\!(1-x^r)^N\omega(t)
T_{\tau_{m-1},n_{m-1}}\!\circ\ldots\circ\, T_{\tau_1,n_1}\!\circ\,
T_{\tau_m,n_m}\left( \frac{x^{\tau_0+\alpha}}{1-(1-x)t}\right)dxdt
$$
with the weight function $\omega(t)$ defined in {\rm(\ref{eq5})}.
\end{lemma}

{\bf Proof.} Applying the following representation introduced by
Pr\'{e}vost  \cite{pr}:
$$
\frac{1}{1-x}+\frac{1}{\log
x}=\int_0^1\frac{\omega(t)}{1-(1-x)t}\,dt,
$$
we have
$$
I(\alpha)=\int_0^1\int_0^1\frac{x^{\tau_0+\alpha}\omega(t)}{1-(1-x)t}A(x)\,dtdx.
$$
As it easily follows  the polynomial $A(x)$ can be written in the
form
$$
A(x)=S_{\tau_1,n_1}\!\circ S_{\tau_2,n_2}\!\circ\ldots\circ
S_{\tau_m,n_m}\left((1-x^r)^N\right).
$$
Since $A(x)$ is symmetric in pairs $(\tau_j,n_j)$ and does not
depend on the order of differential operators $S_{\tau_j,n_j},$ it
is convenient for the sequel to write it as
$$
A(x)=S_{\tau_m,n_m}\!\circ S_{\tau_1,n_1}\!\circ \ldots\circ
S_{\tau_{m-1},n_{m-1}}\left((1-x^r)^N\right).
$$
 Now by Fubini's
theorem and Lemma \ref{l3}, we get the desired equality. \quad\qed

We need also the following simple lemma, which will be used for
estimation purposes.
\begin{lemma} \label{lmax}
Let $a,b,c,d,r,s\in {\mathbb R},$ $r,s,d>0,$ and $b+d\ge a+c\ge
b\ge 0.$ Then the function
$$
f(x,t)=\frac{x^{a+c}(1-x^r)^{sc}t^{c+a-b}(1-t)^{b+d-c-a}}{(1-t+xt)^d}
$$
attains its maximum in $[0,1]\times [0,1]$ at the unique point
$$
x_0=\left(\frac{b}{b+scr}\right)^{\frac{1}{r}}, \qquad t_0=
\frac{c+a-b}{c+a-b+x_0(b+d-a-c)}
$$
and
$$
\underset{0\le x,t\le 1}{\max} f(x,t)=f(x_0,t_0)=
\frac{b^{\frac{b}{r}}(scr)^{sc}(c+a-b)^{c+a-b}(b+d-a-c)^{b+d-a-c}}%
{d^d(b+scr)^{sc+\frac{b}{r}}}.
$$
\end{lemma}

\section{Proof of Theorem \ref{t1}}

\begin{lemma} \label{l4.5}
Let $x,t\in (0,1),$ $\tau_0, n_m, \tau_m\in {\mathbb N}_0,$ and
$\tau_m\le \tau_0\le n_m+\tau_m.$ Then
$$
T_{\tau_m,n_m}\left(\frac{x^{\tau_0}}{1-(1-x)t}\right)=(-1)^{n_m}
\frac{x^{n_m+\tau_m}t^{n_m+\tau_m-\tau_0}(t-1)^{\tau_0-\tau_m}}%
{(1-(1-x)t)^{n_m+1}}.
$$
\end{lemma}

{\bf Proof.} Clearly,
$$
T_{\tau_m,n_m}\left(\frac{x^{\tau_0}}{1-t+xt}\right)=
\frac{x^{n_m+\tau_m}}{n_m!}\left(\frac{x^{\tau_0-\tau_m}}{1-t+xt}\right)^{(n_m)}.
$$
Decomposing the fraction $\frac{x^{\tau_0-\tau_m}}{1-t+xt}$  into
the sum
$$
\frac{x^{\tau_0-\tau_m}}{1-t+xt}=p(x)+\left(\frac{t-1}{t}\right)^{\tau_0-\tau_m}
\frac{1}{1-t+xt},
$$
where $p(x)$ is a polynomial of degree not exceeding
$\tau_0-\tau_m-1,$ and differentiating it $n_m$ times, we get the
required statement. \qed

\begin{lemma} \label{l5}
Under the hypothesis of Theorem {\rm\ref{t1}} one has
\begin{equation}
\begin{split}
&T_{\tau_{m-1},n_{m-1}}\!\circ\ldots\circ\,
T_{\tau_1,n_1}\!\circ\, T_{\tau_m,n_m}\left(
\frac{x^{\tau_0}}{1-(1-x)t}\right) =(-1)^{n_m} \\
&\times\prod_{j=1}^m\binom{n_m+\tau_m-\tau_j}{n_j}
\frac{x^{n_m+\tau_m}t^{n_m+\tau_m-\tau_0}(t-1)^{\tau_0-\tau_m}}%
{(1-t+xt)^{n_m+1}}\,Q_m\left( \frac{xt}{1-t+xt}\right),
\label{eq100}
\end{split}
\end{equation}
where the polynomial $Q_m(y)$ is defined in {\rm(\ref{q})}.
\end{lemma}

{\bf Proof.} If $m=1,$ then (\ref{eq100}) easily follows by Lemma
\ref{l4.5}.
 Suppose $m\ge 2.$ Then consecutive
calculation of the $n_j$th derivatives with respect to $x$  by
Leibniz' rule for $j=1,2,\ldots, m-1$
\begin{equation*}
\begin{split}
&\frac{x^{\tau_j+n_j}}{n_j!}\left(
\frac{t^kx^{n_m+\tau_m+k-\tau_j}}{(1-t+xt)^{n_m+1+k}}\right)^{(n_j)}
=\binom{n_m+\tau_m-\tau_j}{n_j}\frac{x^{n_m+\tau_m}}{(1-t+xt)^{n_m+1}}\\
&\times\sum_{k_j=0}^{n_j}\frac{(-n_j)_{k_j}(n_m+1)_{k+k_j}(1+n_m+\tau_m-\tau_j)_k}%
{\,k_j!\,(n_m+1)_k(1+n_m+\tau_m-\tau_j-n_j)_{k+k_j}}\left(
\frac{xt}{1-t+xt}\right)^{k+k_j}
\end{split}
\end{equation*}
readily leads to the formula (\ref{eq100}). \qed

Now Theorem \ref{t1} easily follows from Lemmas \ref{l4},
\ref{l5}.

\section{Proof of  Theorem \ref{t2}}

If we put $m=1, n_1=cn, \tau_1=an, \tau_0=bn, n\in {\mathbb N},$
in Theorem \ref{t1}, we get
\begin{equation*}
\begin{split}
&\left|\gamma-\frac{1}{r^{cn}}
\sum_{k=0}^{cn}(-1)^{k+cn}\binom{cn}{k}
\binom{rk+(a+c)n}{cn}S_{rk+bn}\right|\\
&\le\frac{1}{r^{cn}}\int_0^1\int_0^1\frac{(1-x^r)^{cn}x^{(a+c)n}t^{(c+a-b)n}
(1-t)^{(b-a)n}\omega(t)}{(1-t+xt)^{cn+1}}\,dxdt \\
&\le\frac{1}{r^{cn}}\left( \underset{0\le x,t\le
1}{\max}f(x,t)\right)^n\int_0^1\int_0^1
\frac{\omega(t)}{1-t+xt}dtdx=\frac{\gamma}{r^{cn}}\left(
\underset{0\le x,t\le 1}{\max}f(x,t)\right)^n
\end{split}
\end{equation*}
with
$$
f(x,t)=\frac{x^{a+c}(1-x^r)^ct^{c+a-b}(1-t)^{b-a}}{(1-t+xt)^c}.
$$
Here we used the fact (see \cite[formula 2.6]{pr}) that
$$
\gamma=\int_0^1\left(\frac{1}{\log x}+\frac{1}{1-x}\right)\,dx.
$$
Now, since $\gamma<1,$ by Lemma \ref{lmax} with $s=1, d=c,$  the
theorem follows. \qed

\section{Proofs of Theorems \ref{t3}, \ref{t4}}

To estimate the speed of convergence of quantities (\ref{3.5}) to
$\gamma$ as $N\to\infty$ we need an upper bound for the polynomial
$Q_m(y).$ In some situations it is possible to get suitable
estimations.

First, we consider the case  $m=2,$ $n_1=n_2,$ $\tau_1=\tau_2.$
Then by Theorem \ref{t1}, we get
\begin{equation*}
\begin{split}
I:&=\left|\frac{(2n_1)!\,r^{2n_1}}{(n_1!)^2}\,\gamma-\sum_{k=0}^{2n_1}(-1)^k
\binom{2n_1}{k}\binom{rk+n_1+\tau_1}{n_1}^2S_{rk+\tau_0}\right| \\
&=\int_0^1\int_0^1\frac{x^{n_1+\tau_1}(1-x^r)^{2n_1}t^{n_1+\tau_1-\tau_0}
(t-1)^{\tau_0-\tau_1}\omega(t)}{(1-t+xt)^{n_1+1}}\,|Q_2(y)|\,dxdt
\end{split}
\end{equation*}
with $y=xt/(1-t+xt).$ The polynomial
$$
Q_2(y)={}\sb 2F\sb
 {1}\left(\left.\atop{-n_1,
n_1+1}{1}\right|y\right)=\frac{1}{n_1!}\left(\frac{d}{dy}\right)^{n_1}
\Bigl(y^{n_1}(1-y)^{n_1}\Bigr)
$$
is a shifted Legendre polynomial $P_{n_1}(u)$ formally identified
as follows:
$$
Q_2(y)=P_{n_1}(1-2y).
$$
By the well-known inequality (see \cite[p.162]{sz})
$$
|P_{n_1}(u)|\le 1, \qquad -1\le u\le 1,
$$
   it follows  that
$$
I\le
\int_0^1\int_0^1\frac{x^{n_1+\tau_1}(1-x^r)^{2n_1}t^{n_1+\tau_1-\tau_0}
(1-t)^{\tau_0-\tau_1}\omega(t)}{(1-t+xt)^{n_1+1}}\,dxdt.
$$
Now, setting $n_1=cn,$ $\tau_1=an,$ $\tau_0=bn$ with $c,b\in
{\mathbb N},$ $a\in {\mathbb N}_0,$ and $0\le b-a\le c,$ we get
\begin{equation*}
\left|\frac{(2cn)!\,r^{2cn}}{((cn)!)^2}\gamma-
\sum_{k=0}^{2cn}(-1)^{k}\binom{2cn}{k}
\binom{rk+(a+c)n}{cn}^2S_{rk+bn}\right| \le\gamma\left(
\underset{0\le x,t\le 1}{\max} f(x,t)\right)^n,
\end{equation*}
where
$$
f(x,t)=\frac{x^{c+a}(1-x^r)^{2c}t^{a+c-b}(1-t)^{b-a}}{(1-t+xt)^c}.
$$
By Lemma \ref{lmax}, the function $f(x,t)$ takes its maximum in
$[0,1]\times[0,1]$ at the unique point $(x_0,t_0),$
at which
$$
f(x_0,t_0)=\frac{b^{\frac{b}{r}}(4cr^2)^c(c+a-b)^{c+a-b}(b-a)^{b-a}}%
{(b+2cr)^{2c+\frac{b}{r}}}.
$$
Since for any positive integer $n\ge 2$
$$
\gamma\frac{(n!)^2}{(2n)!}\le\frac{n}{4^n},
$$
Theorem \ref{t3} follows.  \qed


Another interesting case is described by the following lemma.

\begin{lemma} \label{l6}
Let $n_1,\ldots, n_m\in {\mathbb N},$ $\tau_0, \tau_1,\ldots,
\tau_m\in {\mathbb N}_0,$
 and
$n_m+\tau_m\ge\tau_{j+1}> n_j+\tau_j,$ $j=1,\ldots, m-1.$
  Then
\begin{equation}
\begin{split}
&Q_m(y)=\prod_{j=1}^{m-1}\frac{(n_m+\tau_m-n_j-\tau_j)!}%
{(n_m+\tau_m-\tau_{j+1})!(\tau_{j+1}-n_j-\tau_j-1)!}\,\times \\
&\int_0^1
\!\!\!\!\!\ldots\!\!\int_0^1\prod_{j=1}^{m-1}(1-yu_j\ldots
u_{m-1})^{n_j}u_j^{n_m+\tau_m-\tau_{j+1}}(1-u_j)^{\tau_{j+1}-n_j-\tau_j-1}
du_1\ldots du_{m-1}. \label{eq101}
\end{split}
\end{equation}
Moreover,  $0\le Q_m(y)\le 1$ for $y\in [0,1].$
\end{lemma}

{\bf Proof.} Denoting the integral on the right-hand side of
(\ref{eq101}) by $J$ and substituting
\begin{equation*}
\prod_{j=1}^{m-1}(1-yu_ju_{j+1}\ldots
u_{m-1})^{n_j}=\sum_{k_1=0}^{n_1}\ldots\sum_{k_{m-1}=0}^{n_{m-1}}\prod_{j=1}^{m-1}
\frac{(-n_j)_{k_j}y^{k_j}u_j^{k_1+\ldots+k_j}}{k_j!},
\end{equation*}
we get
\begin{equation*}
\begin{split}
J&=\sum_{k_1=0}^{n_1}\ldots\sum_{k_{m-1}=0}^{n_{m-1}}\prod_{j=1}^{m-1}
\frac{(-n_j)_{k_j}y^{k_j}}{k_j!}\int_0^1u_j^{k_1+\ldots+k_j+n_m+\tau_m-\tau_{j+1}}\times
\\
&\times
(1-u_j)^{\tau_{j+1}-n_j-\tau_j-1}\,du_j=\sum_{k_1=0}^{n_1}\ldots\sum_{k_{m-1}=0}^{n_{m-1}}\prod_{j=1}^{m-1}
\frac{(-n_j)_{k_j}y^{k_j}}{k_j!}\times \\
&\times\frac{\Gamma(k_1+\ldots+k_j+n_m+\tau_m+1-\tau_{j+1})
\Gamma(\tau_{j+1}-n_j-\tau_j)}{\Gamma(k_1+\ldots+k_j+n_m+\tau_m+1-n_j-\tau_j)}
\\
&=\prod_{j=1}^{m-1}\frac{\Gamma(1+n_m+\tau_m-\tau_{j+1})
\Gamma(\tau_{j+1}-n_j-\tau_j)}{\Gamma(1+n_m+\tau_m-n_j-\tau_j)} \\
&\sum_{k_1=0}^{n_1}\ldots\sum_{k_{m-1}=0}^{n_{m-1}}
\prod_{j=1}^{m-1}\frac{(-n_j)_{k_j}(1+n_m+\tau_m-\tau_{j+1})_{k_1+
\ldots+k_j}}{k_j!(1+n_m+\tau_m-n_j-\tau_j)_{k_1+\ldots+k_j}}y^{k_j}
\\
&=\prod_{j=1}^{m-1}\frac{(n_m+\tau_m-\tau_{j+1})!(\tau_{j+1}-n_j-\tau_j-1)!}%
{(n_m+\tau_m-n_j-\tau_j)!}\, Q_m(y).
\end{split}
\end{equation*}
The inequality $0\le Q_m(y)\le 1$ for $y\in [0,1]$ easily follows
from the integral representation (\ref{eq101}).  \qed

Now, Theorem \ref{t4} is a consequence of Theorem \ref{t1} and
Lemma \ref{l6}.

\section{Proof of Theorem \ref{t5}}

Setting $n_j=c_jn,$ $j=1,\ldots,m,$ $C=\sum_{j=1}^mc_j,$
$\tau_1=an+1,$ $\tau_0=bn+m$ in Corollary \ref{c3} we get that the
absolute value of the remainder is less than
$$
\frac{M(\overline{c})}{r^{cn}}\int_0^1\int_0^1\frac{x^{(C+a)n+m}
(1-x^r)^{Cn}t^{(C+a-b)n}(1-t)^{(b+c_m-C-a)n}\omega(t)}%
{(1-t+xt)^{c_mn+1}}\,dxdt
$$
with some constant $M(\overline{c})<C^{m-1},$ since
$$
\prod_{j=1}^{m-1}\frac{Cn+j}{(c_{j+1}+\ldots+c_m)n+m-j}<C^{m-1}.
$$
Denoting
$$
f(x,t)=\frac{x^{C+a}(1-x^r)^Ct^{C+a-b}(1-t)^{b+c_m-C-a}}{(1-t+xt)^{c_m}}
$$
and applying Lemma \ref{lmax} with $s=1, d=c_m,$ we conclude the
theorem. \qed

\end{document}